\documentclass[11pt]{amsart}

\usepackage{amssymb,amscd,amsthm,amsxtra} 
\usepackage{latexsym}
\usepackage{epsfig}

\def\be#1{\begin{equation}\label{#1}}

\newtheorem{thm}{Theorem}[section]

\newtheorem{lem}[thm]{Lemma}

\theoremstyle{definition}
\newcommand{\comment}[1]{}

\newtheorem{defn}[thm]{Definition}
\theoremstyle{remark}
\newtheorem{rem}[thm]{Remark}
\numberwithin{equation}{section}

\begin{document}

\title[Dyadic model of turbulence] {An inviscid dyadic model of turbulence: 
the fixed point and Onsager's conjecture}

\author{Alexey Cheskidov}
\address{Department of Mathematics,
University of Michigan, Ann Arbor, MI  48109-1043}
\email{\tt acheskid@umich.edu}

\author{Susan Friedlander}
\address{Department of Mathematics, Statistics, and
Computer Science, Chicago, IL 60607-7045}
\email{\tt susan@math.uic.edu}

\author{Nata\v{s}a Pavlovi\'{c}}
\address{Department of Mathematics, 
Princeton University, Princeton, NJ 08544-1000}
\email{\tt natasa@math.princeton.edu}

\date{October 3, 2006}

\begin{abstract}
Properties of an infinite system of nonlinearly coupled
ordinary differential equations are discussed. This
system models some properties present in the equations 
of motion for an inviscid fluid such as the skew symmetry
and the 3-dimensional scaling of the quadratic nonlinearity.
It is proved that  the system with forcing has a unique equilibrium
and that {\em every} solution blows up in finite time in $H^{5/6}$-norm.
Onsager's conjecture is confirmed for the model system.
\end{abstract}

\maketitle

\setcounter{equation}{0}
\section{Introduction} 

One of the outstanding questions in fluid dynamics is existence,
uniqueness and regularity of solutions to the Cauchy problem for the three-dimensional 
Euler equations:
\begin{align}
\begin{split} \label{E} 
& \frac{\partial u}{\partial t} = - (u \cdot \nabla) u - \nabla p, \\
& \nabla \cdot u = 0. 
\end{split}
\end{align} 
Important features of the nonlinear term $(u \cdot \nabla) u$ are its bilinearity 
and skew-symmetry. The second property implies conservation of energy for 
sufficiently regular solution of the 
Euler equation (and decay of energy in the context of the Navier-Stokes equations). 
In the past few decades discretized ``toy'' models that preserve the energy properties
of the fluid equations have been proposed and studied both by mathematicians and 
physicists. These models belong to a general class of ``shell'' models which simulate
the energy cascade in turbulent flow. In all of these models the nonlinearity 
of the Euler equations is much simplified by considering only local neigbouring 
interactions between certain scales. However simplifications vary and as a consequence 
the models differ in the number of conserved quantities and in the presence of 
a certain ``monotonicity'' property that we will discuss later. Among  the first 
example of such discretized models is the one introduced by Gledzer \cite{Gledzer} 
in 1970 which was later generalized by  Ohkitani and Yamada \cite{OY}
and is now known as the GOY model. A survey of mathematical developments in 
connection with shell models can be found in the recent book of Bohr et al 
\cite{BJPV}. 

In this article we analyze one of these discretized models, namely 
\begin{equation} \label{Bsystem}
\begin{split} 
\frac{d a_j}{d t} & = 2^{\frac{5}{2}(j-1)} a_{j-1}^{2} - 2^{\frac{5}{2}j} a_{j} a_{j+1} +f_j, \; \; j >0 \\
\frac{d a_0}{d t} & = - a_0 a_1 + f_0.  
\end{split}
\end{equation}
The motivation for this model, and specifically for the scaling appropriate
for 3 dimensions is given in \cite{FP} and \cite{KP} in the case  
$-\infty < j < \infty$ and $f=0$. In \cite{FP} a shell model is 
presented, following the Fourier space analysis of Dinaburg and Sinai 
\cite{DS}, for the 3-dimensional incompressible Euler equations \eqref{E} 
and it is shown that for very specific initial data the vector model reduces
to the scalar system analogous to \eqref{Bsystem} without forcing. 
Furthermore, the wavelet decomposition utilised in \cite{KP} and \cite{KPalpha} 
motivates a discretized model for the 3-dimensional Euler 
equations analogous to \eqref{Bsystem} without forcing.

A derivation of the system ~\eqref{Bsystem} is presented in Section~\ref{derivation}.
The coefficient $a_j^2(t)$ is the {\em total} energy in the frequency space shell
$2^j \leq |k| <2^{j+1}$.
 In this context $l^2$ and $H^s$, respectively energy and
Sobolev norms, are defined as: 
\[ 
\|a(t)\|_{l^2} = \left( \sum_{j = 0}^{\infty} a_j^2(t) \right)^{1/2}, \qquad
\|a(t)\|_{H^{s}} = \left( \sum_{j = 0}^{\infty} 2^{2sj} a_j^2(t) \right)^{1/2}. 
\]
Local in time existence of solutions in $H^{5/2}$ was obtained in \cite{FP} 
using the fixed point techniques that produce local existence in $H^{5/2}$ 
for mild solutions of actual 3 dimensional Euler equations \eqref{E} (see 
Cannone \cite{Can}). Finite time blow-up of \eqref{Bsystem}, without forcing, 
in an appropriate $H^s$ was initially proved in \cite{KP} by exploiting
conservation of energy  and a monotonicity property present in the quadratic
term $a_{j-1}^2$. Variants of the model have been studied in \cite{KZ}, \cite{W},
and \cite{ZH}. 
In particular, it was observed in \cite{W} that \eqref{Bsystem} with $2^{5j/2}$
replaced by $2^j$ could be derived from a discretized version of the
Fourier transform of the 1-dimensional Burger equations under the (unrealistic) restriction 
that there is no spreading of the support.


Variants of the model including dissipation have been studied  \cite{C}, \cite{CLT},
\cite{KP}, \cite{KPalpha} and finite time blow up is proved for certain 
``small'' amounts of dissipation.  More precisely, in \cite{C} a finite time
blow up was proved in the cases where the dissipation degree is such that
the model enjoys the same estimate on the nonlinear term as the Navier-Stokes
equations in the dimensions larger than 4.

Even though a global existence proof from \cite{C} applied to \eqref{Bsystem}
implies that for any initial data in $l^2$ there exists a global in time classical
solution to \eqref{Bsystem}, not all the classical solutions satisfy
the energy balance equation. For instance, regular solutions, i.e., solutions
with bounded $H^{5/6}$-norm, satisfy the energy equality, but the unique
fixed point does not. 

In this present paper we study the system \eqref{Bsystem} with forcing as a
model for 3-dimensional turbulence. The main results are the following:

\begin{itemize}

\item[(a)] For any force $f=(f_0,0,\dots)$, $f_0>0$, there exists a unique fixed point
of the system \eqref{Bsystem}. The fixed point is not in $H^{5/6}$ and is given by
\begin{equation*}
\{a_j\}=\{2^{-5j/6+5/12}\sqrt{f_0}\}.
\end{equation*}

\item[(b)] The system linearized about the fixed point has exponentially
decaying eigenvalues.

\item[(c)]
Every regular solution approaches the fixed point in $l^2$ norm.

\item[(d)] Every solution blows up in finite time in $H^{5/6}$ norm.
\end{itemize}

In a more technical companion paper \cite{CFPatr} we study the 
solutions of \eqref{Bsystem} {\em after} the time of blow-up in 
$H^{5/6}$. Among the properties we prove that 
the $H^s$ norms for $s<5/6$ are locally square integrable in time. 
Moreover, we study the long time behavior of \eqref{Bsystem} 
and prove that the fixed point is a global attractor. This is a
consequence of anomalous or turbulent dissipation that was conjectured
by Onsager \cite{O} and related to Kolmogorov's prediction that in 
fully developed turbulent flow the energy spectrum in the inertial 
range is given by a power law 
\begin{equation} \label{kol} 
E(|k|) = c_0 \bar{\epsilon}^{2/3} |k|^{-5/3}.
\end{equation}

In Section~\ref{Onsager} of the present paper we observe that the energy 
spectrum of the fixed point of the forced model \eqref{Bsystem} 
reproduces Kolmogorov's law \eqref{kol}. Furthermore, regular solutions,
i.e., solutions with bounded
$H^{5/6}$ norm, satisfy the energy equality, whereas the fixed point
is not regular and does not satisfy energy equality. Moreover,
anomalous dissipation cuts in for all the other
solutions as well when the $H^{5/6}$ norm becomes unbounded.



\subsection*{Organization of the paper}

In section 2 we present a derivation of the dyadic model
based on Littlewood-Paley operators. 
In section 3 we review existence and finite time blow-up results
for the model obtained in the past few years. In sections 4 and 5
we prove the results stated above.
In section 6 we discuss how the dyadic model with the 3-dimensional
scaling satisfies Onsager's conjecture and Kolmogorov's 5/3 law.

\subsection*{Acknowledgements} The authors would like to thank 
Marie Farge, Jonathan Mattingly, Kai Schneider and Eric Vanden-Eijnden
for very helpful discussions. S.F. was partially supported by NSF grant 
number DMS 0503768. N.P. was partially supported by NSF grant 
number DMS 0304594.

\setcounter{equation}{0}
\section{Model based on Littlewood-Paley decomposition}
\label{derivation}
Let us start by considering 3D Euler equations with 
zero force (for the sake of simplicity):
\be{Euler}
{\partial u \over \partial t} + (u \cdot \nabla) u + \nabla p
= 0.
\end{equation}

\subsection{Euler equations on $\mathbb{R}^3$.} \label{s:deriv}
Given a velocity $u(x,t)$, let
\[
\hat{u}(\xi,t)=\int_{\mathbb{R}^3} e^{-ix\cdot\xi} u(x,t) \, dx, \qquad \xi \in \mathbb{R}^3,
\]
be its Fourier transform.
We shall use a standard Littlewood Paley decomposition. 
More precisely, we consider Fourier multipliers
$P_j$ (on $L^2({\mathbb R}^3)$)   
$$ \widehat{P_j u}(\xi,t) = p_j(\xi) \hat u(\xi,t),$$
so that their symbols $p_j(\xi)$ are smooth and supported
in ${2 \over 3} 2^j < |\xi| < 3 (2^j)$ and such that 
$p_j(\xi)=p_0(2^{-j} \xi)$ and $\sum_j p_j(\xi)=1.$
Also we denote by $\tilde P_j$ the following multiplier
$\tilde P_j = \sum_{k=-2}^{2} P_{j+k}$. 

We apply the Leray projection $T$ onto the divergence free vectors followed 
by an application of $P_j$ on the equation \eqref{Euler} to obtain: 
\be{PjEuler} 
{\partial P_j u \over \partial t} + P_j(T(u \cdot \nabla u)) = 0. 
\end{equation}
Our goal is to see what we can say about $\|P_j u\|_{L^2}$
because that is the relevant quantity in the definition of 
the Sobolev norm 
\be{sob-act} 
\|f\|_{H^s} = \left( \sum_{j=-\infty}^\infty (1+ 2^{2js}) \|P_j f\|_{L^2}^2 \right)^{1/2}.  
\end{equation} 
We use Bony's paraproduct formula for 
the nonlinear term. More precisely, 

$$P_j(T(u \cdot \nabla u)) = N_{j,lh} + N_{j,hl} + N_{j,hh} + N_{loc},$$
where the low-high part is given by
$$N_{j,lh}=\sum_{j'< j - 4} P_j((P_{j'} u) \cdot \tilde P_j \nabla u),$$
the high-low part is given by
$$N_{j,hl}=\sum_{j'< j - 4} P_j((\tilde P_j u )\cdot P_{j'} \nabla u),$$
the high-high part is given
$$N_{j,hh}=
\sum_{j'> j + 4}
P_j((\tilde P_{j'} u) \cdot P_{j'} \nabla u)
+ \sum_{j'> j + 4}
P_j(( P_{j'} u) \cdot \tilde P_{j'} \nabla u),$$
and the local part is given via
$$N_{loc}=
\sum_{ j-4 \leq j' \leq  j+4} 
P_j((\tilde P_{j'} u) \cdot P_{j'} \nabla u)
+ \sum_{ j- 4 \leq j' \leq j+4} 
P_j(( P_{j'} u) \cdot \tilde P_{j'} \nabla u).$$

Now let us start modeling. We recall that the Fourier multiplier
$P_jf(x)$ are given via 
$$ \widehat{P_j f}(\xi) = p_j(\xi) \hat{f}(\xi)$$ 
has the symbol $p_j(\xi)$ which is, roughly speaking, supported 
for $|\xi| \sim 2^j$. This combined with
the belief that only local frequency scales are relevant 
in a turbulence cascade, motivates us to keep only local 
interactions in the nonlinear term, i.e. we keep only the
equivalent of a modified version of $N_{loc}$ in the model that 
we propose below. More precisely, we introduce a model 
that describes the evolution of the coefficients 
\be{coef} 
a_j(t) = \|P_j u(x,t)\|_{L^2}. 
\end{equation} 
Here $a_j^2(t)$ represents the total energy in the shell.
In order to introduce in the model an analogue of the $N_{loc}$ expression,
we need to see what scaling comes out of 
$ \|P_j u \cdot \nabla P_{j'}(u)\|_{L^2}$ when $j$ and $j'$ are close. 
By H\"{o}lder's inequality combined with the fact that 
$p_j(\xi)$ are supported for $|\xi| \sim 2^j$ we obtain: 

\begin{align} 
\|P_ju \cdot \nabla P_{j'} u \|_{L^2} 
& \leq \|P_j u \|_{L^\infty} \|\nabla P_{j'} u\|_{L^2} \\
& \leq 2^{\frac{3}{2}j} \|P_j u \|_{L^2}\cdot 2^{j'} \|P_{j'} u \|_{L^2} \label{ber}\\
& = 2^{\frac{5}{2}j}a_ja_{j'}, \label{scal}  
\end{align} 
where to obtain \eqref{ber} we use Bernstein's inequality which can be stated on ${\mathbb{R}}^n$
in terms of Littlewood-Paley operators as follows  
$$ \|P_j u\|_{L^q} \leq 2^{(\frac{1}{p} - \frac{1}{q})nj} \|P_j u\|_{L^p}, \mbox{ for } q > p.$$ 

Now we are ready to propose the following model for the 3D Euler equations:  
\be{dyadicR3} 
\frac{da_j}{dt} = 2^{\frac{5}{2}(j-1)} a_{j-1}^2 - 2^{\frac{5}{2}j} a_{j} a_{j+1},
\qquad j \in \mathbb{Z}.
\end{equation}  
Note that due to \eqref{coef}, the $L^2({\mathbb{R}}^3)$ norm of the
fluid velocity $u$ is
\[
\|u\|_{L^2}= \left( \sum_{j=-\infty}^\infty |a_j|^2 \right)^{1/2}.
\] 
The Sobolev norm \eqref{sob-act} of $u$ is 
\be{sob}
\|u\|_{H^s} = \left( \sum_{j=-\infty}^\infty (1+ 2^{2js}) |a_j|^2 \right)^{1/2}.  
\end{equation} 

\subsection{Space periodic 3D Euler equations.}
Consider now Euler equations in a periodic box $\Omega=[0,L]^3$. By Galilean
change of variables we can assume that the space average of $u$ is zero.
Given a velocity $u(x,t)$, let
\[
\hat{u}_k(t)=\frac{1}{L^3}\int_{\Omega} e^{-i\frac{2\pi}{L}x\cdot k} u(x,t)\, dx, \qquad k \in \mathbb{Z}^3,
\]
be its Fourier coefficients. Define $a_j^2(t)$ to be the total energy in the shell
$2^j\leq|k|<2^{j+1}$
\begin{equation} \label{ajFcoef}
a_j(t)^2= \sum_{2^{j}\leq |k| < 2^{j+1} } |\hat{u}_k(t)|^2, \qquad j =0,1,2,\dots.
\end{equation}
An analysis similar to that given in Subsection~\ref{s:deriv} motivates the following model of the 3D Euler equations:
\begin{equation}
\begin{split}
\label{dyadic} 
\frac{da_j}{dt} &= 2^{\frac{5}{2}(j-1)} a_{j-1}^2 - 2^{\frac{5}{2}j} a_{j} a_{j+1},
\qquad j =1,2, \dots,\\
\frac{d a_0}{d t} & = - a_0 a_1.
\end{split}
\end{equation}  
Due to \eqref{ajFcoef}, the $L^2(\Omega)^3$ and $H^s(\Omega)^3$ norms of the fluid
velocity $u$ can be defined as
\begin{equation} \label{norms-periodic}
\|u\|_{L^2}= \|a\|_{l^2}=\left( \sum_{j=0}^\infty |a_j|^2\right)^{1/2}, \qquad
\|u\|_{H^s} = \left( \sum_{j=0}^\infty 2^{2js} |a_j|^2 \right)^{1/2}. 
\end{equation}

\section{Global existence and blow-up} \label{glob} 

We write the dyadic model \eqref{dyadic} as
\begin{align} 
\begin{split} \label{dsystem}
\frac{d a_j}{d t} & = \lambda^{j-1} a_{j-1}^{2} - \lambda^j a_{j} a_{j+1} +f_j, \; \; j >0,\\
\frac{d a_0}{d t} & = - a_0 a_1+f_0,
\end{split} 
\end{align} 
with the forcing $f_j\geq 0$ for all $j$.

For notational simplicity we adopt $\lambda^j$ as the scaling
parameter in the computations in sections \ref{glob} and \ref{uniq}. We do
this to illustrate that the results are qualitively independent of the
exact choice of lambda ( which depends on the spatial dimension
and the construction of the model). As we discussed in sections 1
and 2 the relevant lambda for our 3 dimensional dyadic model
is $2^{5/2}$. This exponent determines the values of the critical exponents in
the Sobolev space results proved below.


We say that $a(t)=(a_0(t),a_1(t),\dots)$ is a solution to \eqref{dsystem} if it is
a classical solution in the usual sense. More precisely, we have the following.
\begin{defn}
A solution on $[0,T]$ (or $[0, \infty)$, if
$T=\infty$) of \eqref{dsystem} is an $l^2$-valued
function $a(t)$ defined for $t \in [0, T]$, such that $a_j \in C^1([0,T])$
and $a_j(t)$ satisfies \eqref{dsystem} on $(0,T)$ for all $j$.
\end{defn}
Note that if $a(t)$ is a solution on $[0,T]$, then automatically
$a_j \in C^{\infty}([0,T])$ for all $j$. The global existence of a classical solution for
any initial data $a(0) \in l^2$ was proved in Cheskidov \cite{C} for a closely related
system. The idea of the proof is as follows, see Theorem 4.1 in \cite{C}
for more details. Given an arbitrary time interval $[0,T]$,
consider a sequence of Galerkin approximations to \eqref{dsystem}.
It is easy to show that this sequence is weakly equicontinuos on $[0,T]$.
Therefore, thanks to Ascoli-Arzela theorem, we can pass to a convergent
subsequence and obtain the existence of a weak
solution to \eqref{dsystem} on $[0,T]$. Since the nonlinear term has a finite number of terms, this weak solution is also a classical solution. Hence we have the following
result.
\begin{thm}
For every $a^0 \in l^2$ and $f\in l^2$, there exists a solution $a(t)$ to
\eqref{dsystem} on $[0,\infty)$ with $a(0) = a^0$.
\end{thm}

Define the $H^s$ norm of $a$ as
\[  \|a(t)\|_{H^{s}} = \left( \sum_{j = 0}^{\infty} 2^{2sj} a_j^2(t) \right)^{1/2}.\]
Due to \eqref{norms-periodic}, $\|\cdot\|_{H^s}$ can be viewed as the model
analogue of the Sobolev norm of a fluid velocity.
For each scale $j$ we define the energy as 
\be{jen} E_{j}(t) = \frac{1}{2}a_{j}^{2}(t),
\end{equation}
and the energy of the whole system 
\be{en} E(t) = \frac{1}{2} \|a(t)\|^2_{l^2}=\sum_{j=0}^{\infty} E_{j}(t). 
\end{equation}
We define the energy of a ``box'' $B_{J}$ as
\be{boxen} E_{B_{J}}(t)= \sum_{j=J}^{\infty} E_{j}(t).
\end{equation}

Formally multiplying \eqref{dsystem} by $a_j$ and summing over all $j$ gives 
\be{ensum} 
\frac{1}{2}\frac{d}{dt} \sum_{j=0}^{\infty} a_j^2 = 
\sum_{j = 1}^{\infty} \lambda^{j-1} a_{j-1}^2 a_j - \sum_{j=0}^{\infty} \lambda^j a_j^2 a_{j+1}
+ \sum_{j = 0}^{\infty} f_ja_j. 
\end{equation}  
We recall that $\lambda=2^{5/2}$.
Hence, if $\|a(t)\|_{H^{5/6}}$ is bounded on some  time interval, then the summations
on the right hand side of \eqref{ensum} are uniformly 
convergent, and we obtain the following energy balance property:
\be{consen} 
\frac{d}{dt} E(t) =  \sum_{j = 0}^{\infty} f_ja_j. 
\end{equation}

This motivates the following.
\begin{defn}
A solution $a(t)$ of \eqref{dsystem} is called regular (or strong) on $[T_1,T_2]$ if $\|a(t)\|_{H^{5/6}}$
is bounded on $[T_1,T_2]$.
\end{defn}
Hence, the total energy of a regular solution to \eqref{dsystem} satisfies the
energy equality \eqref{consen}.
A regular solution also satisfies 
\be{monotbox} 
\frac{d}{dt} E_{B(J)}(t) = \lambda^{J-1} a_{J-1}^2 a_J, \; \; J>0.
\end{equation} 

Treating the $j$-th equation in \eqref{dsystem} as an ODE for $a_j(t)$  gives 
\be{integ} 
\begin{split}
a_j(t) &= \frac{1}{\mu(t)} \left[ a_j(0) + \int_{0}^{t} \lambda^{j-1} a_{j-1}^{2}(\tau) \mu(\tau) \; d\tau \right], \; \; j>0,\\
a_0(t) &= \frac{1}{\mu(t)} \left[ a_0(0) + f_0\int_{0}^{t} \mu(\tau) \; d\tau \right],
\end{split}
\end{equation} 
where 
\begin{equation}
\mu(t) = \exp \left(\int_{0}^{t} \lambda^j a_{j+1}(\tau) \; d\tau \right).
\end{equation}
Hence if $a_j(0) \geq 0$ for all $j$, we conclude that this property is 
preserved by evolution of \eqref{dsystem} and 
$$a_j(t) \geq 0 \mbox{ for all } j \geq 0.$$
Thus, \eqref{monotbox} implies that energy in each box 
$B_J$ of a regular solution increases monotonically in time, i.e. 
\be{posenchange} 
\frac{d}{dt} E_{B_J}(t) \geq 0, \; \; J>0. 
\end{equation} 
The inequality \eqref{posenchange} implies that the system produces a
successive cascade  of energy into higher and higher scales $J$.   

The basic properties of system \eqref{dsystem} described above contribute to the proofs of certain results 
that have appeared in Katz-Pavlovi\'{c} \cite{KP}, Friedlander-Pavlovi\'{c} \cite{FP}, Kiselev-Zlato\v{s} \cite{KZ}. 
We recall these theorems.

\begin{thm} \label{locex} 
Let $a^0 \in H^s$ for some $s\geq 5/2$ with $a^0_j \geq 0$
for all $j \geq 0$. Then 
there exists a  time $T = T(\|a^0\|_{H^s}) > 0$ such that a unique solution 
$a(t)$ to \eqref{dsystem} with $a(0)=a^0$ exists and $a(\cdot) \in C([0,T]; H^s)$.
\end{thm} 
This theorem is proved in Friedlander-Pavlovi\'{c} \cite{FP} using Picard's fixed 
point argument in the case $f=0$.
In particular, this theorem implies local existence of regular solutions with
initial data in $H^{5/2}$. The local existence of regular solutions is not yet known
for initial data not in $H^{5/2}$.

\begin{thm} \label{blow-upH} 
Let $a(t)$ be a solution of  \eqref{dsystem} with $a(t) \in H^1$, $a_j(0) \geq 0$ for all $j \geq 0$, and $a(0) \ne 0$. 
Then $\|a(t)\|_{H^s}$ becomes infinite in finite time for all $s>5/6$. 
\end{thm}
The first proof of this type of blow-up result in $H^s$ for a system 
closely related to \eqref{dsystem} with $f=0$ was given in Katz-Pavlovi\'{c}
\cite{KP}. Versions for $H^s$ blow up were given in 
Friedlander-Pavlovi\'{c} \cite{FP}, Kiselev-Zlato\v{s}s \cite{KZ} (the above theorem)
and for a Navier-Stokes model in \cite{KP} and Cheskidov \cite{C}.  

Crucial to these blow-up proofs is the cancellation of the infinite sums
on the right hand side of \eqref{ensum} leading to the energy properties
\eqref{consen} and \eqref{monotbox}. In the context of the scaling of
this present paper, the smallest value of $s$ for which finite time blow-up is proved is
$s=5/6$. We note that the claim by Waleffe \cite{W} to use such a proof
to produce blow-up down to $s>0$ (for a variant of the model with $2^j$ scaling) 
is not justified because the telescoping sums in \eqref{ensum} do not converge for
$0<s<s_0$, where $s_0$ is a relevant critical Sobolev exponent ($s_0=1/3$ for $2^j$
scaling).

Consider now a forcing of the form $f_0>0$ and $f_j=0$ for all $j\geq 1$.
Clearly,
\[
a_j = \lambda^{-\frac{j}{3}+\frac{1}{6}} \sqrt{f_0}
\]
is a fixed point of \eqref{dsystem}. Since $\lambda = 2^{5/2}$, note that the fixed point is not regular, i.e.,
it is not in $H^{5/6}$. Note also that it does not satisfy the energy balance equation
\eqref{consen}. Indeed, for a fixed point we have
\[
0 = \frac{d}{dt} E(t) < \sum_{j = 0}^{\infty} f_ja_j. 
\]
Moreover, at the end of this paper we will prove that there is no global in time
regular solution to \eqref{dsystem}, i.e., every solution $a(t)$ with $a(0) \in l^2$
blows up in finite time in $H^{5/6}$-norm. This result improves
Theorem~\ref{blow-upH}.

\setcounter{equation}{0}
\section{Uniqueness of the fixed point} \label{uniq}

We continue studying the dyadic model
\begin{equation} \label{stabsystem}
\begin{split}
\frac{d a_j}{d t} & = \lambda^{j-1} a_{j-1}^{2} - \lambda^j a_{j} a_{j+1}+f_j , \; \; j >0\\
\frac{d a_0}{d t} & = - a_0 a_1 + f_0.
\end{split}
\end{equation}
For the forcing of the form $f_0>0$ and $f_j=0$ for all $j\geq 1$
we will investigate the spectrum of the system linearized about the
fixed point
\be{fix} 
\{a_j\} = \{\lambda^{-\frac{j}{3}+\frac{1}{6}} \sqrt{f_0}\}.                  
\end{equation}

\subsection{Uniqueness of the fixed point} 
\begin{thm}
For every force $f=(f_0,0,0, \dots)$, $f_0>0$, there exists a
unique fixed point of \eqref{stabsystem}
\begin{equation}
\{a_j\} = \{\lambda^{-\frac{j}{3}+\frac{1}{6}} \sqrt{f_0}\}.
\end{equation}
\end{thm}
\begin{proof}
Consider
\begin{equation}
A_j = \lambda^{\frac{j}{3}-\frac{1}{6}} f_0^{-\frac{1}{2}} a_j.
\end{equation}
Then the equations \eqref{stabsystem} for a fixed point become
\begin{align*}
A_{j-1}^2 - A_jA_{j+1} &=0, \; \; j>0,\\
A_0A_1&=1.
\end{align*}
Clearly,
\begin{equation}
A_j = 1, \qquad j \geq 0
\end{equation}
is a solution to the above system, which corresponds to
\begin{equation}
a_j = \lambda^{-\frac{j}{3}+\frac{1}{6}} \sqrt{f_0}.
\end{equation}
Let us show that this is the only fixed point. Indeed, if $\{a_j\}$ is
a fixed point, then
\begin{equation}
A_0A_1=1, \qquad \frac{A_{j-1}}{A_j} = \frac{A_{j+1}}{A_{j-1}}, \qquad j \geq 1.
\end{equation}
From these equalities we obtain
\begin{equation*}
\frac{A_0}{A_1}=\frac{A_2}{A_0}, \qquad
\frac{A_0}{A_2}=\frac{A_2A_3}{A_0A_1}, \qquad
\frac{A_0}{A_3}=\frac{A_3A_4}{A_0A_1}, \qquad
\frac{A_0}{A_4}=\frac{A_4A_5}{A_0A_1}, \qquad \dots.
\end{equation*}
Thus, we have
\begin{equation}
\frac{A_0}{A_j} =A_jA_{j+1}, \qquad j\geq 2.
\end{equation}
Therefore,
\begin{equation}
A_{j+1} = A_0 A_j^{-2}, \qquad j \geq 2.
\end{equation}
Hence,
\begin{equation}
A_{j}=\left\{
\begin{split}
&A_0^{(1+2^{j-2})/3}  A_2^{-2^{(j-2)}}, &j\geq 3 \mbox{ odd},\\
&A_0^{(1-2^{j-2})/3}  A_2^{2^{(j-2)}},  &j\geq 3 \mbox{ even}.
\end{split}
\right.
\end{equation}
Finally, since $A_2= A_0^3$, we obtain
\begin{equation}
A_{j}=\left\{
\begin{split}
&A_0^{-3\cdot2^{(j-2)}+(1+2^{j-2})/3}, &j\geq 3 \mbox{ odd},\\
&A_0^{3\cdot 2^{(j-2)}+(1-2^{j-2})/3},  &j\geq 3 \mbox{ even}.
\end{split}
\right.
\end{equation}
Clearly, $A_0$ has to be equal to $1$ in order for $\{a_j\}$ to be bounded. 
Then $A_j=1$ for all $j$.
\end{proof}

Note that the above argument can also be applied to a more general force
of the form $f=(f_0,f_1, \dots, f_k,0,0,\dots)$ with $f_k>0$.


\begin{lem}
Let $f=(f_0,f_1, \dots, f_k,0,0,\dots)$ with $f_k>0$.
If there exists a fixed point $\{a_j\}$ of \eqref{stabsystem}, then
\begin{equation}
a_j=\lambda^{\frac{1}{3}} a_{j+1}, \qquad \forall j \geq k.
\end{equation}
\end{lem}

\begin{proof}
Let $a_j$ be a fixed point. Since $f_k>0$, we have that
$a_ka_{k+1} >0$. Consider
\begin{equation}
A_j = \frac{\lambda^{(j-k)/3-1/6} a_j}{\sqrt{a_ka_{k+1}}}.
\end{equation}
Then $A_k A_{k+1}=1$ and
\begin{equation}
A_jA_{j+1}=A_{j-1}^2, \qquad j\geq k+1.
\end{equation} 
Therefore, we have
\begin{equation}
\frac{A_{j-1}}{A_j} = \frac{A_{j+1}}{A_{j-1}}, \qquad j \geq k+1.
\end{equation}
From these equalities we obtain
\begin{equation*}
\frac{A_k}{A_{k+1}}=\frac{A_{k+2}}{A_k}, \quad
\frac{A_k}{A_{k+2}}=\frac{A_{k+2}A_{k+3}}{A_{k}A_{k+1}}, \quad
\frac{A_k}{A_{k+3}}=\frac{A_{k+3}A_{k+4}}{A_kA_{k+1}}, \quad
\frac{A_k}{A_{k+4}}=\frac{A_{k+4}A_{k+5}}{A_kA_{k+1}}, \ \dots
\end{equation*}
Hence, it follows that
\begin{equation}
\frac{A_k}{A_j} =A_j A_{j+1}, \qquad j\geq k+2.
\end{equation}
Therefore,
\begin{equation}
A_{j+1} = A_k A_j^{-2}, \qquad j \geq k+2.
\end{equation}
Thus,
\begin{equation}
A_{j}=\left\{
\begin{split}
&A_k^{(1+2^{j-2})/3} A_{k+2}^{-2^{(j-2)}}, &j\geq k+3, \ j-k \mbox{ odd},\\
&A_k^{(1-2^{j-2})/3}  A_{k+2}^{2^{(j-2)}},  &j\geq k+3, \ j-k \mbox{ even}.
\end{split}
\right.
\end{equation}
Finally, since $A_{k+2}= A_k^3$, we obtain
\begin{equation}
A_{j}=\left\{
\begin{split}
&A_k^{-3\cdot 2^{(j-2)}+(1+ 2^{j-2})/3}, &j\geq k+3, \ j-k \mbox{ even},\\
&A_k^{3\cdot 2^{(j-2)}+(1-2^{j-2})/3},  &j\geq k+3, \ j-k \mbox{ odd}.
\end{split}
\right.
\end{equation}
Clearly, $A_k$ has to be equal to $1$ in order for $\{a_j\}$ to be bounded. 
Then $A_j=1$ for all $j\geq k$, which implies that $a_j=\lambda^{1/3}a_{j+1}$
for all $j\geq k$.
\end{proof}

Using this lemma one can prove the following theorem.

\begin{thm}
Let $f=(f_0,f_1, \dots, f_k,0,0,\dots)$ with
$f_j\geq0$ for all $j$. Then there exists a
unique fixed point $\{a_j\}$ of \eqref{stabsystem}. Moreover,
\begin{equation}
a_j= \lambda^{-\frac{j}{3}}C, \qquad \forall j \geq k,
\end{equation}
for some nonnegative constant $C$.
\end{thm}

\subsection{Spectral stability} 
We consider the fixed point of \eqref{consen} with $f=(f_0,0,0,\dots)$, $f_0>0$.
Rescaling the variables, we can assume that $f_0=\lambda^{-1/3}$ for simplicity.
Then the fixed point is of the form $\{\lambda^{-j/3}\}$. To investigate its spectral stability we write
\be{pert} 
a_j(t) = \lambda^{-\frac{j}{3}} + \epsilon b_j.
\end{equation} 
where the perturbation $\{b_j(t)\}$ satisfies
\begin{align} 
\frac{d b_j}{dt} & = \lambda^{\frac{2j}{3}} 
( 2\lambda^{-\frac{2}{3}} b_{j-1} - \lambda^{-\frac{1}{3}} b_j - b_{j+1} ) 
+ \epsilon (\lambda^{j-1} b_{j-1}^2 - \lambda^j b_j b_{j+1})\; \; j>0 \label{bj} \\
\frac{d b_0}{dt} & = -\lambda^{-\frac{1}{3}} b_0 - b_1 - \epsilon b_0 b_1. \label{b0}
\end{align} 
The linearized equations for small $\epsilon$ are
\begin{align}  
\frac{d b_j}{dt} & = \lambda^{\frac{2j}{3}} 
( 2\lambda^{-\frac{2}{3}} b_{j-1} - \lambda^{-\frac{1}{3}} b_j - b_{j+1} )\; \; j>0 \label{lbj} \\
\frac{d b_0}{dt} & =  -\lambda^{-\frac{1}{3}} b_0 - b_1. \label{lb0}
\end{align} 

We seek a solution to \eqref{lbj}-\eqref{lb0} of the form 
\be{lam} 
b_j = c_j e^{\mu t}
\end{equation} 
with $\mu$ real.  

\begin{lem} \label{lem-lambda}
There are no positive eigenvalues to the system \eqref{lbj} - \eqref{lb0}. Moreover,
there exist solutions in $H^s$, $s < 5/6$ to the system \eqref{lbj} - \eqref{lb0}
of the form \eqref{lam} with the eigenvalues $\mu$ negative. 
\end{lem}

\begin{proof} 
Substitution of \eqref{lam} into \eqref{lbj}-\eqref{lb0} gives 
\begin{align} 
& \lambda^{-\frac{1}{6}} c_{j+1} + \alpha_j c_j - 2\lambda^{-\frac{5}{6}} c_{j-1} = 0, \; \; j>0 \label{cj} \\ 
& \alpha_0 c_0 + \lambda^{-\frac{1}{6}} c_1 = 0, \label{c0} 
\end{align} 
where 
\be{al} 
\alpha_j = \lambda^{-\frac{1}{6}} (\lambda^{-\frac{1}{3}} + \lambda^{-\frac{2j}{3}} \mu), \; \; j \geq 0.
\end{equation} 
We shall construct a sequence $\{c_j\}$, $c_j \neq 0$, $j \geq 0$  that solves the system \eqref{cj}-\eqref{c0}.

Define 
\be{dj} 
d_j = \lambda^{-\frac{1}{6}} \frac{c_{j+1}}{c_j}, \; \; j \geq 0.
\end{equation}
It follows from \eqref{cj} and \eqref{dj} that
\be{recdj}
d_{j-1} = \frac{1}{2^{3/2}(\alpha_j + d_j) } \; \; j \geq 1.
\end{equation} 
Hence 
\be{dn} 
d_n = [\alpha_{n+1}, \alpha_{n+2}, ...],
\end{equation} 
where $[\alpha_{n+1}, \alpha_{n+2}, ...]$ denotes the 
continued fraction 
$$ [\alpha_{n+1}, \alpha_{n+2}, ...] 
= \cfrac{ 1 }{ 2^{3/2}\alpha_{n+1} + \cfrac{ 2^{3/2} }
                                   { 2^{3/2}\alpha_{n+2} + \cfrac{2^{3/2}}{...}    
                                   }
             }.$$

Observe that from \eqref{al} 
\be{lim-al} 
\lim_{j \rightarrow \infty} \alpha_j = \lambda^{-\frac{1}{2}},
\end{equation}
for any finite $\mu$. 
Clearly for sufficiently large $n$ the coefficients $\alpha_{n+m}$, $m \geq 1$  
in \eqref{dn} are positive and the infinite continued fraction \eqref{dn} 
is convergent to a positive number. Furthermore it follows that 
\be{lim-dn} 
d_{\infty} = \lim_{n \rightarrow \infty} d_n = \lambda^{-\frac{1}{2}} < 1.
\end{equation}

From \eqref{cj} - \eqref{c0} we obtain the characteristic equation 
\be{char} 
-\alpha_0 = [\alpha_1, \alpha_2, ...],
\end{equation} 
where $\alpha_j$'s are given by \eqref{al}. We observe that \eqref{char} 
has no roots $\mu \geq 0$, since in this case all $\alpha_j$ including 
$\alpha_0$ are positive. 

Suppose there exists $\mu$ real and negative such that 
\eqref{char} is satisfied. Choose
\begin{align} 
\begin{split} \label{cs} 
c_0 & = 1 \\
c_1 & = \lambda^{\frac{1}{6}} d_0 \\
...\\
c_j & = \lambda^{\frac{j}{6}} d_{j-1} d_{j-2} ...d_0.
\end{split} 
\end{align} 
Then this sequence satisfies \eqref{cj}-\eqref{c0} by construction
and in view of \eqref{lim-dn} 
\be{asym-cj} 
c_j \sim \lambda^{-\frac{j}{3}}, \mbox{ for } j >>1. 
\end{equation}
Thus the $l^2$ norm of $\{c_j\}$ is finite. In fact
the $H^s$ norm of $\{c_j\}$ is finite for $s < 5/6$. 

We are now going to study the characteristic equation 
\eqref{char}. Let $X(\mu)$ denote the function 
\be{X} 
X(\mu) = \alpha_0(\mu) + [\alpha_{1}(\mu), \alpha_{2}(\mu), ...], 
\end{equation}
with $\mu \in (-\infty, 0)$. 
Then to solve \eqref{char} we look for zeros of the function $X(\mu)$. 

Since $\alpha_j(\mu)$ is given by \eqref{al} we can write $X(\mu)$
in the form 
\be{recX} 
X(\mu) = 
\lambda^{-\frac{1}{6}} (\lambda^{-\frac{1}{3}} + \mu) + \frac{1}{2^{3/2} X( \lambda^{-2/3} \mu ) }.
\end{equation} 

\begin{figure}
\centering
\includegraphics[width=3in]{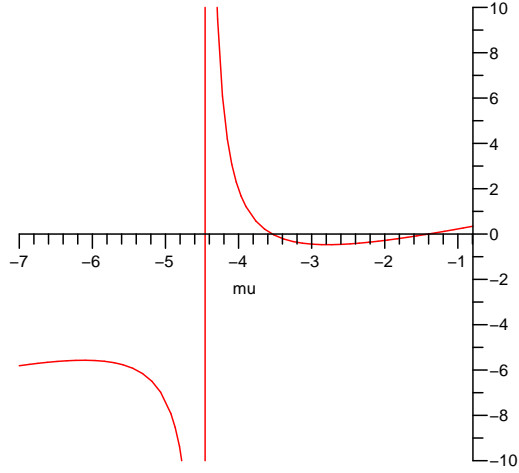}
\caption{Graph of $X(\mu)$}
\label{fig:xgraph2006Sep19.eps}
\end{figure} 
From \eqref{recX} and \eqref{al} we observe that 
\be{X0}
X(0) > 0,
\end{equation} 
and 
\be{Xinf} 
 \lim_{\mu \rightarrow -\infty} X(\mu) = -\infty.
\end{equation}
It follows from \eqref{X0} and \eqref{Xinf} that $X(\mu)$ cannot be
continuous and positive for all $\mu \in(-\infty, 0)$. For $\mu \in (-\lambda^{-1/3},0)$
each entry of the continuous fraction defining $X(\mu)$ in \eqref{char} is positive.
Hence $X(\mu)$ is continuous in $(-\lambda^{-1/3},0)$. Let $(\mu_0, 0)$ be the
maximal interval of continuity of $X(\mu)$. Then $X(\mu)$ is
continuous at $\lambda^{-2/3}\mu_0$. This implies that $X(\lambda^{-2/3}\mu_0)=0.$
Figure 1 gives a plot using Maple to graph the function $X(\mu)$ and show
that the zeros of $X(\mu)$ are approximately $\mu = -1.4$ and 
$\mu = -3.6$. 
\end{proof}

\section{Behaviour of regular solutions} \label{bofregsol}

In this section we show that every regular solution (i.e. with bounded
$H^{5/6}$ norm) approaches the fixed point in the $l^2$ norm linearly
in time. Furthermore, we prove that every solution with initial datum
in $l^2$ blows up in finite time in the $H^{5/6}$ norm.

\begin{thm}\label{stabilitylemma}
Let $a(t)$ be a regular solution of \eqref{stabsystem} on $[T_1,T_2]$
with $f=(\lambda^{-1/3},0,0,\dots )$. We write
$
a_j(t) = \lambda^{-\frac{j}{3}} + b_j(t).
$
Then
\[
\|b(T_2)\|_{l^2} - \|b(T_1)\|_{l^2} \leq -\frac{1}{2}\cdot\lambda^{-\frac{1}{3}}(T_2-T_1).
\]
\end{thm}

\begin{proof}
Note that $b(t)$ satisfies the following system of equations:
\begin{equation}
\begin{split} \label{rpbj}
\frac{d b_j}{dt} & = \lambda^{\frac{2j}{3}} 
( 2\lambda^{-\frac{2}{3}} b_{j-1} - \lambda^{-\frac{1}{3}} b_j - b_{j+1})
 +  ( \lambda^{j-1} b_{j-1}^2 - \lambda^j b_j b_{j+1} ), \; \; j\geq 1, \\
\frac{d b_0}{dt} & = -\lambda^{-\frac{1}{3}} b_0 - b_1 -  b_0 b_1.
\end{split} 
\end{equation}
Multiplying it by $b_j$ and taking a sum from $j=0$ to $j=k$ we obtain
\begin{multline} \label{beftel-enbj}
 \frac{1}{2} \frac{d}{dt} \sum_{j=0}^{k} b_j^2 
= -\lambda^{-\frac{1}{3}} b_0^2 - b_0b_1 - b_0^2 b_1 \\
 + \sum_{j=1}^{k} \left( 
\lambda^{\frac{2j}{3}} ( 2\lambda^{-\frac{2}{3}} b_{j-1}b_j - \lambda^{-\frac{1}{3}}b_j^2 - b_jb_{j+1}) 
+   ( \lambda^{j-1} b_{j-1}^2b_j - \lambda^j b_j^2 b_{j+1} ) \right). 
\end{multline}
It now follows that
\begin{multline} \label{enbj}
\frac{1}{2} \frac{d}{dt} \sum_{j=0}^{k} b_{j}^{2}   
=  -\left(\lambda^{-\frac{1}{3}} b_0^2 + b_0b_1\right) \\
 + \sum_{j=1}^{k} \lambda^{\frac{2j}{3}} 
\left( 2\lambda^{-\frac{2}{3}} b_{j-1}b_j - \lambda^{-\frac{1}{3}}b_j^2 - b_jb_{j+1}\right)
- \lambda^k b_k^2b_{k+1}.
\end{multline} 
Also we can rewrite \eqref{enbj} as 
\begin{multline}\label{enbjrew} 
\frac{1}{2} \frac{d}{dt} \sum_{j=0}^{k} b_{j}^{2} 
= - \lambda^{-\frac{1}{3}} \sum_{j=0}^{k} \lambda^{\frac{2j}{3}} b_j^2 + 
\sum_{j=0}^{k} \lambda^{\frac{2j}{3}} b_j b_{j+1}\\
-2\lambda^{\frac{2k}{3}} b_kb_{k+1}
 - \lambda^k b_k^2b_{k+1}. 
\end{multline} 
However,
\[
\begin{split} 
\sum_{j=0}^{k} \lambda^{\frac{2j}{3}} 
(\lambda^{-\frac{1}{6}} b_j - \lambda^{\frac{1}{6}} b_{j+1})^2
&= \sum_{j=0}^{k} \lambda^{\frac{2j}{3}} 
(\lambda^{-\frac{1}{3}} b_j^2 + \lambda^{\frac{1}{3}} b_{j+1}^2 - 2 b_j b_{j+1}) \\
& = -2 \left[ - \lambda^{-\frac{1}{3}} \sum_{j=0}^{k} \lambda^{\frac{2j}{3}} b_j^2   
+ \sum_{j=0}^{k} \lambda^{\frac{2j}{3}} b_j b_{j+1} \right]\\
& \quad - \lambda^{-\frac{1}{3}} b_0^2
+\lambda^{\frac{2k}{3}+\frac{1}{3}}b_{k+1}^2.  
\end{split}
\] 
Hence \eqref{enbjrew} gives 
\begin{multline} \label{ebjrew1}
\frac{1}{2} \frac{d}{dt} \sum_{j=0}^{k} b_{j}^{2} 
 = -\frac{1}{2} \left[
\sum_{j=0}^{k} \lambda^{\frac{2j}{3}} 
(\lambda^{-\frac{1}{6}} b_j - \lambda^{\frac{1}{6}} b_{j+1})^2 
+ \lambda^{-\frac{1}{3}} b_0^2
\right]\\
 + \frac{1}{2}\lambda^{\frac{2k}{3} + \frac{1}{3}}b_{k+1}^2
-2\lambda^{\frac{2k}{3}} b_kb_{k+1}
 - \lambda^k b_k^2b_{k+1}.
\end{multline}

Now note that since $a(t)$ is regular on $[T_1,T_2]$, we have that
\begin{equation}
\lim_{j \to \infty} \int_{T_1}^{T_2} \lambda^{\frac{2j}{3}} a^2_j(t) \, dt =0.
\end{equation}
Hence, since $b_j=a_j-\lambda^{-j/3}$, we obtain
\begin{multline} \label{ebjrew1-1}
\lim_{j \to \infty} \int_{T_1}^{T_2} \left( \frac{1}{2}\lambda^{\frac{2k}{3}+\frac{1}{3}}b_{k+1}^2
-2\lambda^{\frac{2k}{3}} b_kb_{k+1} - \lambda^k b_k^2b_{k+1}\right) \, dt\\
\begin{split}
&=(T_2-T_1)\left(\frac{1}{2} \lambda^{-\frac{1}{3}} - 2\lambda^{-\frac{1}{3}} + \lambda^{-\frac{1}{3}} \right)\\
&=-\frac{1}{2}\cdot\lambda^{-\frac{1}{3}}(T_2-T_1).
 \end{split}
\end{multline}
Now \eqref{ebjrew1} and \eqref{ebjrew1-1} imply that
\begin{equation}
\|b(T_2)\|_{l^2} - \|b(T_1)\|_{l^2} \leq -\frac{1}{2}\cdot\lambda^{-\frac{1}{3}}(T_2-T_1).
\end{equation}
\end{proof}


\subsection{Blow-up in finite time}

Using the results of the previous subsection, we prove that there
is no global in time regular solution to \eqref{stabsystem} with $f=(\lambda^{-1/3},0,0,\dots)$.

\begin{thm} \label{t:Blow-up}
All the solution to \eqref{stabsystem} with $f=(\lambda^{-1/3},0,0,\dots)$
blow up in finite time in $H^{5/6}$-norm. More precisely, for every solution $a(t)$
with $a(0) \in l^2$, there exists a time
\[
0\leq t^*\leq 2\lambda^{\frac{1}{3}} \|b(0)\|_{l^2},
\]
such that
\[
\limsup_{t \to t^*} \|a(t)\|_{H^{5/6}} = \infty.
\]
\end{thm}
\begin{proof}
First note that every $a_j(t)$ is continuous, i.e., $a(t)$ is weakly in $l^2$
continuous. Therefore, if $\|a(t^*)\|_{H^{5/6}}=\infty$, then
\begin{equation}
\limsup_{t \to t^*-} \|a(t)\|_{H^{5/6}} = \limsup_{t \to t^*+} \|a(t)\|_{H^{5/6}}=\infty.
\end{equation}

Now assume that the statement of the theorem is not true. Then there exists
a solution $a(t)$ to \eqref{stabsystem}
which is regular on $[0,T]$ with $T= 2\lambda^{1/3} \|b(0)\|_{l^2}$. Then Theorem~\ref{stabilitylemma}
implies that
\begin{equation}
\begin{split}
\|b(T)\|_{l^2} &\leq  \|b(0)\|_{l^2} -2^{-1}\lambda^{-\frac{1}{3}}T\\
&=0.
\end{split}
\end{equation}
Therefore,
\begin{equation}
a_j(T)=\lambda^{-\frac{j}{3}},
\end{equation}
i.e., $\|a(T)\|_{H^{5/6}}=\infty$. Due to the remark in the beginning of the proof, 
\begin{equation}
\limsup_{t\to T-} \|a(t)\|_{H^{5/6}} = \infty,
\end{equation}
a contradiction.
\end{proof}

Note that there is no global in time regular solution to \eqref{stabsystem} with a more general force $f=(f_0,0,0,\dots)$, $f_0>0$. This follows by application of
Theorem~\ref{t:Blow-up} to a rescaled system \eqref{stabsystem}.

\section{Onsager's conjecture and Kolmogorov's $5/3$rd law.}
\label{Onsager}
In the past 60 years since its original presentation much attention has
been given to the statistical theories of turbulence developed by
Kolmogorov and Onsager. For recent reviews see, for example Eyink and
Sreenivasan \cite{ES}, Robert \cite{R}. Onsager \cite{O} conjectured a clear
distinction between 2 and 3 dimensions: in 2 dimensions the energy of a
turbulent flow is conserved, however, in 3 dimensions dissipation of
energy persists in the limit of vanishing viscosity. This phenomenon is
now referred to as turbulent or anomalous dissipation. It is suggested
that an appropriate mathematical description of 3-dimensional turbulent
flow is given by weak solutions of the Euler equations which are not
regular enough to conserve energy. Kolmogorov's theory \cite{K} predicts
that in a fully developed turbulent flow the energy spectrum $E(|k|)$ in
the inertial range is given by
\be{Kolm}
E(|k|)=c_0\bar{\epsilon}^{2/3} |k|^{-5/3},
\end{equation}
where $\bar{\epsilon}$ is the average of the energy
dissipation rate. This law stated in physical space means that the Holder exponent $h$
of the velocity is $1/3$ in a statistically averaged sense. Onsager conjectured that
for exponents $h>1/3$ the energy is conserved and that this ceases to be true for
$h\leq 1/3$. An elegant proof of the conservation of energy of weak solutions
of the 3 dimensional Euler equations in Besov spaces $B_3^{h,\infty}$, $h>1/3$,
was given by Constantin et al \cite{CET}.

The model system that we study in this present paper reproduces the phenomenon
described above. Regular solutions, i.e., solutions with bounded $H^{5/6}$ norm,
satisfy the energy equality, whereas the fixed point
\begin{equation} \label{O-fpoint}
\{a_j\}=\{2^{-5j/6+5/12}\sqrt{f_0}\}.
\end{equation}
does not.
Note that a necessary condition for a function in physical space to have Holder
exponent $h$ is that the spectral exponent is greater or equal to $1+2h$
(i.e., the spectral exponent $5/3$ of the fixed point in our model is exactly
the Onsager critical exponent corresponding to $h=1/3$).

The model~\eqref{Bsystem} is derived under the assumption that $a_j^2(t)$ is the
{\em total} energy in the frequency space shell $2^j \leq |k| <2^{j+1}$. Hence,
by \eqref{O-fpoint}, the energy spectrum $E(|k|)$ for the fixed point is given by
the expression
\begin{equation}
E(|k|)=2^{5/6}f_0|k|^{-5/3}.
\end{equation}
Since the dissipation rate for the fixed point is equal to the energy
input rate, we have
\[
\epsilon=a_0f_0= 2^{5/12}f_0^{3/2}.
\]
Thus the fixed point satisfies Kolmogorov's law \eqref{Kolm} with
$c_0=2^{5/9}$.

In Section~\ref{bofregsol} we proved that every regular solution approaches
the fixed point and, furthermore, blows up in
finite time in $H^{5/6}$ norm. After this time (as we show in \cite{CFPatr})
the solutions are not regular enough to conserve energy and anomalous
dissipation for an inviscid system produces Kolmogorov's turbulent
energy spectrum.

Anomalous dissipation and ``life after blow-up" has been examined recently
in a {\em linear} 1-dimensional model where the system is "simple" enough to
be fully solvable but still produces anomalous dissipation via energy
cascade towards higher wave numbers \cite{MS}. In our more technical
companion paper \cite{CFPatr} we examine ``life after blow-up" for the
system~\eqref{Bsystem}. In particular, we prove that the fixed point is
a global attractor. Since the support of any time-average measure belongs
to the global attractor, the average dissipation rate is equal to the dissipation
rate of the fixed point. Thus Kolmogorov's law stated above for the fixed
point is valid for the system.

\begin{rem}
It is well known that the 1-dimensional inviscid Burgers equation produces energy
dissipation for weak solutions with shocks (as we noted, a rescaled version
of \eqref{Bsystem} can be motivated by Burgers equation in Fourier space).
In contrast, the Cauchy problem for the 3-dimensional incompressible Euler
equations remains open (and a major challenge). The 3-dimensional
Euler equations may (or may not) share with the model system \eqref{Bsystem}
the turbulent processes conjectured by Onsager for 3-dimensional fluids.
The model {\em per se} could be written with n-dimensional scaling and the system
itself has no behavior that is inherently 3-dimensional or incompressible.
We note, however, that the reduction of an incompressible 3-dimensional
vector model for the fluid equations to the scalar system \eqref{Bsystem}
given in \cite{FP} cannot be implemented in two dimensions. Furthermore,
any such model that does not preserve enstrophy cannot be a suitable model for the
2-dimensional Euler equations.
\end{rem}


\begin{thebibliography}{99}


\bibitem{BJPV} T. Bohr, M. Jensen, G. Paladin and A. Vulpiani:
Dynamical Systems Approach to Turbulence,
Cambridge University Press, 1998.

\bibitem{Can} M. Cannone:
Harmonic analysis tools for solving the incompressible Navier-Stokes equations,
{\em Handbook of mathematical fluid dynamics} Vol. III, 161--244, 
North-Holland, Amsterdam, 2004. 

\bibitem{C} A. Cheskidov, 
 Blow-up in finite time for the dyadic model of the Navier-Stokes equations,
{\em Transactions of AMS}, to appear.  


\bibitem{CF} A. Cheskidov and C. Foias, On global attractors of the 3D
Navier-Stokes equations, {\it J. Diff. eq.}, to appear.

\bibitem{CFPatr} A. Cheskidov, S. Friedlander, N. Pavlovi\'{c},
An inviscid dyadic model of turbulence: the global attractor, 
{\em In preparation}.

\bibitem{CET} P. Constantin, W. E, and E. S. Titi, OnsagerÕs conjecture on the energy conservation for solutions of EulerÕs equation, {\it Comm. Math. Phys.} {\bf{165}}
(1994), 207-209.

\bibitem{CLT} P. Constantin, B. Levant and E. Titi,
Analytic study of the shell model of turbulence, 
{\em Preprint}. 

\bibitem{DS} E.I. Dinaburg and Ya. G. Sinai,
A quasilinear  approximation for the three-dimensional Navier-Stokes system,
{\em Moscow Mathematical Journal} {\bf{1}}, No. 3 (2001), 381--388.

\bibitem{ES} G. L. Eyink and K. R. Sreenivasan,
Onsager and the theory of hydrodynamic turbulence, {\it Rev. Mod. Phys.} {\bf{78}}
(2006).

\bibitem{FP} S. Friedlander and N. Pavlovi\a'{c},
Blow up in a three dimensional vector model for the Euler equations,
{\em Commun. Pure Appl. Math.} {\bf{57}}, No.6 (2004), 705--725.

\bibitem{Gledzer} E. B. Gledzer,
System of hydrodynamic type admitting two quadratic integrals of motion, 
{\em Sov. Phys. Dokl.} {\bf{18}}, No. 4 (1973), 216--217.

\bibitem{KP} N. Katz and N. Pavlovi\a'{c},
Finite time blow-up for a dyadic model of the Euler equations,
{\em Transactions of AMS} {\bf{357}}, No. 2 (2005), 695-708. 

\bibitem{KPalpha} N.H. Katz and N. Pavlovi\'{c},
A cheap Caffarelli-Kohn-Nirenberg inequality for
the Navier Stokes equation with hyper-dissipation,
{\em GAFA} {\bf{12}} (2002), 355--379.

\bibitem{KZ} A. Kiselev and A. Zlato\v{s},
On discrete models of the Euler equation, 
{\em IMRN} {\bf{38}} (2005) No. 38, 2315--2339.

\bibitem{K} A. N. Kolmogorov, The local structure of tur- 
bulence in incompressible viscous fluids at very large 
Reynolds numbers, {\it Dokl. Akad. Nauk. SSSR} {\bf{30}} (1941), 301--305.

\bibitem{MS} J. C. Mattingly, T. Suidan, and E. Vanden-Eijnden,
Simple systems with anomalous dissipation and energy cascade,
{\it preprint}.

\bibitem{OY} K. Ohkitani and M. Yamada,
Temporal intermittency in the energy cascade process and 
local Lyapunov analysis in fully developed model turbulence,
{\em Prog. Theor. Phys.} {\bf{81}}, No. 2 (1989), 329--341.

\bibitem{O} L. Onsager, Statistical Hydrodynamics, {Nuovo Cimento (Supplemento)} {\bf{6}} (1949),  279--287.

\bibitem{R} R. Robert, Statistical Hydrodynamics ( Onsager revisited ),
{\it Handbook of Mathematical Fluid Dynamics}, vol 2 ( 2003),
1--55. Ed. Friedlander and Serre. Elsevier.

\bibitem{W} F. Waleffe
On some dyadic models of the Euler equations, 
{\em Proc. Amer. Math. Soc.} {\bf{134}} (2006), 2913--2922. 

\bibitem{ZH} V. Zimin and F. Hussain, Wavelet based model for small-scale turbulence,
{\it Phys. Fluids} {\bf{7}} (1995), 2925--2927.


\end{thebibliography}
\end{document}